%Format: amsppt
\input amstex
\documentstyle{amsppt}

\input label.def
\input degt.def
\input dd.def

\input epsf
\def\picture#1{\epsffile{#1-bb.eps}}

%\input debug.def
%\PrintLabels

{\catcode`\@=11
\gdef\proclaimfont@{\sl}}

\def\dash{\item"\hfill--\hfill"}
\def\Dashes{\widestnumber\item{--}}

\def\ie,{\emph{i.e.},}
\def\eg,{\emph{e.g.},}
\def\etc{\emph{etc}}
\def\cf.{\emph{cf\.}}
\def\via{\emph{via}}

\def\Maple{{\tt Maple}}
\def\GAP{{\tt GAP}}
\def\Term#1{$\DG{#1}$}
\def\term{\Term{10}}

\loadbold

\def\bA{\bold A}
\def\bD{\bold D}
\def\bE{\bold E}
\def\bJ{\bold J}
\def\bU{\bold U}

\def\bY{\bold Y}
\def\bW{\bold W}

\let\Ga\alpha
\let\Gb\beta
\let\Gg\gamma
\let\Gd\delta
\let\Ge\epsilon
\let\Gs\sigma

\def\BB{\bar B} % the trigonal curve
\def\LL{\bar L} % the section
\def\PP{\bar P} % image of a singular point $P$
\let\GG\pi % the group of the curve
\def\tX{\tilde X} % the resolved covering $K3$-surface
\def\tc{\tilde c} % the involution on $\tX$
\def\DG#1{\Bbb D_{#1}} % the dihedral group
\def\CG#1{\Z_{#1}} % the cyclic group
\def\QG{Q_8} % the quaternion group
\def\GAPcolon{\rtimes}
\def\discr{\operatorname{discr}}
\def\Pic{\operatorname{Pic}}
\def\CK{\Cal K} % the kernel of the expension $S\subset\tilde S$
\def\Cp#1{\Bbb P^{#1}}
\def\tSigma{\tilde\Sigma}
\def\tS{\tilde S}
\def\Sk{\operatorname{Sk}}

\def\1{^{-1}}
\def\even{_{\roman{even}}}
\def\odd{_{\roman{odd}}}
\def\ord{\operatorname{ord}}
\def\ls|#1|{\mathopen|#1\mathclose|}

\def\one{\circ}
\def\two{\.\circ}
\def\four{\bullet}

\topmatter

\author
Alex Degtyarev
\endauthor

\title
On irreducible sextics with non-abelian fundamental group
\endtitle

\address
Department of Mathematics,
Bilkent University,
06800 Ankara, Turkey
\endaddress

\email
degt\@fen.bilkent.edu.tr
\endemail

\abstract
We calculate the fundamental groups $\pi=\pi_1(\Cp2\sminus B)$ for
all irreducible plane sextics~$B\subset\Cp2$
with simple singularities for
which $\pi$ is known to admit a dihedral quotient $\DG{10}$.
All groups found
are shown to be finite, two of them being of large order: $960$
and $21600$.
\endabstract



\endtopmatter

\document

\section{Introduction}

\subsection{Motivation}
Recall that a plane sextic $B\subset\Cp2$ is said to be of
\emph{torus type}, or \emph{tame}, if its equation can be
represented in the form $p^3+q^2=0$, where $p$ and~$q$ are some
homogeneous polynomials of degree~$2$ and~$3$, respectively.
Essentially, sextics of torus type were introduced by O.~Zariski
as the ramification loci of cubic surfaces. The fundamental group
$\pi=\pi_1(\Cp2\sminus B)$ of an irreducible
sextic~$B$ of torus type is known to
be infinite; in particular, it is nonabelian. For a long time, no
other examples of nonabelian groups were known, which lead M.~Oka
to a conjecture~\cite{Oka.conjecture} that the fundamental group
of an irreducible sextic that is not of torus type is always abelian.
The conjecture was disproved in~\cite{degt.Oka}, \cite{degt.Oka2},
and for the counterexamples \emph{for which the group~$\pi$ was
computed} it turned out to be finite, with the exception of one
family with non-simple singularities. Besides, it was also shown
that, for each irreducible sextic that is not of torus type, the
abelinization of the commutant of~$\pi$ is finite. (This assertion
is a restatement of the proved part of Oka's conjecture, related
to the Alexander polynomial.) Thus, the following
statement seems to be a reasonable replacement for the original
conjecture.

\theorem[Conjecture]\label{conjecture}
Let $B$ be an irreducible plane sextic with simple singularities
and not of torus type. Then the fundamental group
$\pi_1(\Cp2\sminus B)$ is finite.
\endtheorem

The fundamental groups of all irreducible sextics with a
non-simple singular point are found in~\cite{degt.Oka}
(the case of a quadruple point)
and~\cite{degt.Oka2}
(the case of a singular point adjacent to~$\bJ_{10}$).
On the other hand, the construction
of sextics with simple singularities and nonabelian
fundamental group~$\pi=\pi_1(\Cp2\sminus C)$
suggested in~\cite{degt.Oka} is rather indirect; it proves
that $\pi$ has a dihedral quotient, but it is not suitable to
compute~$\pi$ exactly. In this paper, we attempt to substantiate
Conjecture~\ref{conjecture}
by computing the groups of some of the
curves discovered in~\cite{degt.Oka}.

\subsection{Principal results}
For a group~$G$, denote by $G'=[G,G]$ its commutant, or derived
group, and let $G''=(G')'$ \etc. We use the notation $\CG{n}$
and $\DG{n}$ for, respectively, the cyclic and dihedral groups of
order~$n$.

To shorten the statements, we introduce the term
\emph{generalized \Term{2n}-sextic} to stand for a plane sextic~$B$
whose fundamental group $\pi_1(\Cp2\sminus B)$
%admits a dihedral quotient
factors to~$\DG{2n}$, $n\ge3$.
A \emph{\Term{2n}-sextic} is an
irreducible generalized
\Term{2n}-sextic with simple singularities.
Recall, see~\cite{degt.Oka}, that there are \Term6-, \Term{10}-,
and \Term{14}-sextics; all \Term6-sextics are of torus type, and
the \term-sextics
form eight equisingular
deformation families, one family for each of the following sets of
singularities:
$$
\gather
4\bA_4,\quad 4\bA_4\oplus\bA_1,\quad 4\bA_4\oplus2\bA_1,\quad
4\bA_4\oplus\bA_2,\\
\bA_9\oplus2\bA_4,\quad \bA_9\oplus2\bA_4\oplus\bA_1,\quad
\bA_9\oplus2\bA_4\oplus\bA_2,\quad
2\bA_9.
\endgather
$$
The objective of the paper is the computation of the fundamental
groups of all \term-sextics.
The principal result is the following
theorem.

\theorem\label{th.main}
Let $B$ be a \term-sextic.
Then the group $\GG=\pi_1(\Cp2\sminus B)$ is finite.
Furthermore, one has $\GG=\DG{10}\otimes\CG3$ with the
following two exceptions\rom:
\roster
\item\local{4A4+2A1}
The set of singularities of~$B$ is $4\bA_4\oplus2\bA_1$\rom: then
$\ord\GG=960$ and one has $\GG''/\GG'''=(\CG2)^4$ and $\GG'''=\CG2$.
\item\local{A9+2A4+A2}
The set of singularities of~$B$ is
$\bA_9\oplus2\bA_4\oplus\bA_2$\rom: then $\ord\GG=21600$ and
$\GG''$ is the only perfect group of order~$720$,
see~\cite{HoltPlesken}.
\endroster
\rom(In all cases, including the exceptional ones,
$\GG/\GG'=\CG6$ and $\GG'/\GG''=\CG5$.\rom)
\endtheorem

According to~\cite{degt.Oka}, the \Term{14}-sextics form two
equisingular deformation families, with the sets of singularities
$3\bA_6$ and $3\bA_6\oplus\bA_1$. For the former family, the group
has recently been shown to be $\DG{14}\times\CG3$,
see~\cite{degt-Oka}.
The group of the sextics with the set of singularities
$3\bA_6\oplus\bA_1$ (which are all projectively equivalent) is
still unknown.

The groups of the \term-sextics with the sets of
singularities $4\bA_4$ and $4\bA_4\oplus\bA_1$ were
found independently by C.~Eyrol and M.~Oka~\cite{Oka.D10}.

\subsection{Other results and contents}
The starting point of the computation is
Theorem~\ref{th.reduction}, which provides an explicit geometric
construction for \term-sextics. (In~\cite{degt.Oka}, only the
existence of \term-sextics with the sets of singularities listed
above is proven.) We show that each \term-sextic~$B$ is
a double covering of a very particular rigid trigonal
curve~$\BB$ in the Hirzebruch surface~$\Sigma_2$; the curve~$\BB$
has two type~$\bA_4$ singular points, and various sets of
singularities for~$B$ are obtained by varying the ramification
locus. Theorem~\ref{th.reduction} is dealt with in
Section~\ref{S.construction}.

The fundamental groups are computed in Section~\ref{S.groups}:
we merely apply the classical approach
due to van Kampen to the ruling of the
Hirzebruch surface~$\Sigma_2$. In a few difficult cases, the
resulting representations are studied using \GAP.

Section~\ref{S.equations} is not directly related to
Theorem~\ref{th.main}: we use the representation given by
Theorem~\ref{th.reduction} to produce explicit equations for
\term-sextics. The equation depends on three
parameters $a,b,c\in\C$, see~\eqref{eq.sextic}; we analyze the
parameter space and describe the triples~$(a,b,c)$
resulting in particular sets of singularities.

\subsection{Acknowledgements}
I am thankful to the organizers of the
Fourth Franco-Japanese Symposium on Singularities held in Toyama
in August, 2007 for their hospitality and for the excellent
working conditions.
A great deal of time consuming computations used in the paper
was done using \GAP\ and \Maple, and I am taking this opportunity
to extend my gratitude to the creators of these indispensable software
packages.

\section{The construction}\label{S.construction}

\subsection{The reduction}
Recall that any involutive automorphism $c\:\Cp2\to\Cp2$ has a
fixed line $L=L_c$ and an isolated fixed point $O=O_c$, and the
quotient $\Cp2(O)\!/c$ is the rational geometrically ruled
surface~$\Sigma_2$.
(Here, $\Cp2(O)$ stands for the plane~$\Cp2$ with $O$ blown up.)
The images in~$\Sigma_2$ of~$O$ and~$L$ are,
respectively, the exceptional section~$E$ and a generic
section~$\LL$, so that $\Cp2(O)$ is the double covering
of~$\Sigma_2$ ramified at $\LL+E$.

Recall also that the semigroup of classes of effective divisors
on~$\Sigma_2$ is generated by~$E$ and the class~$F$ of a fiber of
the ruling. One has $E^2=-2$, $F^2=0$, $F\circ E=1$, and
$K_{\Sigma_2}=-2E-4F$.

The principal result of this section is the following theorem.

\theorem\label{th.reduction}
Let $B\subset\Cp2$ be a \term-sextic.
Then $\Cp2$ admits an
involution~$c$ preserving~$B$, so that $O_c\notin B$,
%. The quotient $\Cp2/c$ is a
%quadratic cone, so that its resolution of singularities is the
%rational ruled surface~$\Sigma_2$. The image in~$\Sigma_2$ of the
%exceptional divisor over the fixed point~$O$ is the exceptional
%section $E\subset\Sigma_2$, and the image of the fixed line~$L$ is
%a generic section $\LL\subset\Sigma_2$. The image of the
and the image of~$B$ in $\Sigma_2=\Cp2(O_c)\!/c$
is a trigonal curve~$\BB$, disjoint from~$E$,
with two type~$\bA_4$ singular points.
\endtheorem

Theorem~\ref{th.reduction} is proved at the end of this section,
in~\ref{s.proof} below. Theorem~\ref{th.reduction} has a partial
converse, see Theorem~\ref{th.converse}.

\remark{Remark}
The proof of Theorem~\ref{th.reduction} given in~\ref{s.proof}
uses the theory of $K3$-surfaces; we do
prove that \emph{each} \term-sextic is
symmetric. The fact that each of the eight equisingular deformation
families of such sextics admits a symmetric \emph{representative}
can easily be established using the results of
Section~\ref{S.groups}, where each set of singularities
listed in the introduction is realized by a
symmetric \term-sextic.
\endremark

\subsection{The covering $K3$-surface}\label{s.K3}
Let $B$ be a plane sextic with simple singularities. Consider the
double covering $X\to\Cp2$ ramified at~$B$ and its minimal
resolution~$\tX$. Then, $\tX$ is a $K3$-surface and the deck
translation $\tau\:\tX\to\tX$ of the covering is a holomorphic
anti-symplectic (\ie, reversing holomorphic $2$-forms) involution.

Recall that $H_2(\tX)\cong2\bE_8\oplus3\bU$ is an even unimodular
lattice of rank~$22$ and signature~$-16$.
For a singular point~$P$ of~$B$,
denote by $D_P\subset H_2(\tX)$ the set of classes of
exceptional divisors
%(or $(-2)$-curves)
over~$P$; we use the same
notation~$D_P$ for the incidence graph of these divisors, which
is an irreducible Dynkin diagram of the same name
$\bA$--$\bD$--$\bE$ as the type of~$P$. Note that, for a
type~$\bA$ singular point, the action of~$\tau_*$ on~$D_P$ is the
only nontrivial symmetry of the graph.
Let $\Sigma_P\subset H_2(\tX)$ be the sublattice spanned by~$D_P$;
it is an irreducible negative definite root system.

Denote $\Sigma=\bigoplus\Sigma_P$, the summation running
over all singular
points~$P$ of~$B$, and let $S=\Sigma\oplus\Z h\subset H_2(\tX)$,
where
$h\in H_2(\tX)$ is the class realized by the pull-back of a generic
line in~$\Cp2$. One has $h^2=2$, and the sums above are
orthogonal.
%Denote by $D=\bigcup_PD_P$ the basis of~$\Sigma$
%formed by the exceptional divisors, as well as its intersection
%graph.
Let $\tSigma\subset\tS\subset H_2(\tX)$ be the primitive hulls
of~$\Sigma$ and~$S$, respectively. The finite index extension
$\tS\supset S$ is determined by its \emph{kernel}~$\CK$, which is
an isotropic subgroup of the discriminant group $\discr S$.
(For the definition of the discriminant group and its relation to
lattice extensions, see V.~V.~Nikulin~\cite{Nikulin}.) As shown
in~\cite{JAG}, if $B$ is irreducible, then
$\CK\subset\discr\Sigma$. According to~\cite{degt.Oka}, it
is the kernel~$\CK$ that essentially enumerates the dihedral
quotients of
$\pi_1(\Cp2\sminus B)$.

\subsection{Symmetric sextics}\label{s.symmetric}
Let $B$ be a plane sextic with simple singularities, and denote
by~$\tX$ the covering $K3$-surface. Let $\tc\:\tX\to\tX$ be a
holomorphic symplectic (\ie, preserving holomorphic $2$-forms)
involution. As is known, $\tc$ has eight fixed points, and the
quotient $Y=\tX'\!/\tc$ is again a $K3$-surface, where $\tX'$
is~$\tX$ with the fixed points of~$\tc$ blown up.

Since the projection $\tX\to\Cp2$ is the map
defined by the linear
system $h\in\Pic\tX$, the two involutions~$\tc$, $\tau$ commute if
and only if
the induced automorphism~$\tc_*$ of $H_2(\tX)$
preserves~$h$. In this case, $\tc$ descends to
an involution $c\:\Cp2\to\Cp2$ which
preserves~$B$. Let $O=O_c$ and
$L=L_c$. In what follows,
{\it we always
assume that $B$ does not contain~$L$ as a component.} We
fix the notation~$\LL$ and~$\BB$ for the images
of~$L$ and~$B$, respectively, in~$\Sigma_2$.

Alternatively, if $\tc_*(h)=h$, then $\tau$ descends
to an anti-symplectic involution $\tau_Y\:Y\to Y$, and the
quotient $Y\!/\tau_Y$ blows down to~$\Sigma_2$.

\lemma\label{image.class}
Let $O$, $B$, $\BB$, \etc\. be as above.
If $O\notin B$, then
$\BB\in\ls|3E+6F|$ is a trigonal
curve disjoint from~$E$. If $O\in B$, then
$\BB\in\ls|2E+6F|$ is a
hyperelliptic curve, $\BB\circ E=2$.
\endlemma

\proof
The branch locus of the ramified covering $Y\to\Sigma_2$ consists
of~$\BB$, $\LL$, and, if~$O\in B$, the exceptional section~$E$. On
the other hand, since~$Y$ is a $K3$-surface, the branch locus is
an anti-bicanonical curve, \ie, it belongs to~$\ls|4E+8F|$. Since
$\LL\in\ls|E+2F|$, the statement follows.
\endproof

\lemma\label{image.trivial}
Let $P$ be a $c$-invariant type~$\bA_p$ singular point
of~$B$,
%and denote by~$\PP$ its image in~$\Sigma_2$.
and let $\PP\in\Sigma_2$ be its image.
Assume that $p>1$ and that $\tc_*$ acts
trivially on~$D_P$. Then $P\in L$ and $\PP$
is a type $\bA_{2p+1}$ singular point of $\BB+\LL$,
\ie, a
point of $(p+1)$-fold intersection of~$\LL$
and~$\BB$ at a smooth branch of~$\BB$.
Conversely, the double covering of a curve~$\BB$ as above has
a type~$\bA_p$ singular point.
\endlemma

\proof
Since each curve in~$D_P$ is preserved by~$\tc$ (as a set), the
intersection points of the curves must be fixed by~$\tc$.
Furthermore, each of the two outermost curves must contain one
more fixed point of~$\tc$. Blowing up the fixed points, one
obtains
a sequence of rational curves with the following incidence graph:
$$
\DDcenter
@(\one)@---@(\four)@---@(\one)@---@(\four)\withdots\length7@---
 @(\four)@---@(\one)@---@(\four)@---@(\one)
\endDD
$$
(Here, $\DDcenter@(\one)\endDD$ and
$\DDcenter@(\four)\endDD$
stand, respectively, for $(-1)$- and
$(-4)$-curves; the vertices in the diagram alternate and their
total number is $2p+1$.) The $(-1)$-curves are in the fixed point
set of~$\tc$, and the $(-4)$-curves are not. Hence, the projection
of the exceptional divisor to~$Y$ is a sequence of $(2p+1)$
rational curves whose incidence graph~$D$ is~$\bA_{2p+1}$. The
involution~$\tau_Y$ induced from~$\tau$ acts on~$D$ as the only
nontrivial symmetry; in particular, it is nontrivial on the middle
curve. Thus, the projection to~$Y\!/\tau_Y$ is a chain of $p$
$(-2)$-curves ending in a $(-1)$-curve; it blows down to a single
type~$\bA_{2p+1}$ singular point of the branch locus.

The converse statement can be proved by analyzing a local equation
of~$\BB+\LL$, so that $\LL=\{y=0\}$, and substituting
$y\mapsto y^2$.
\endproof

\lemma\label{image.nontrivial}
Let $P$ be a $c$-invariant type~$\bA_p$ singular point
of~$B$,
%and denote by~$\PP$ its image in~$\Sigma_2$.
and let $\PP\in\Sigma_2$ be its image.
Assume that either $p=1$ or $p>1$ and $\tc_*$ acts nontrivially
on~$D_P$. Then $p=2k-1$ is odd and
either
\roster
\item\local{PinL}
$P\in L$ and $\PP$ is a type~$\bD_{k+2}$
\rom(type~$\bA_3$ if $p=1$\rom) singular point of $\BB+\LL$, or
\item\local{P=O}
$P=O$ and $\PP$ is a type~$\bD_{k+1}$ singular point
\rom(a pair of type~$\bA_1$ singular points if $p=1$\rom)
of $\BB+E$.
\endroster
The type~$\bD_s$ singular point above is formed by the section~$\LL$
or~$E$ intersecting~$\BB$ with multiplicity~$2$ at a
type~$\bA_{s-3}$ singular point of~$\BB$.
Conversely, the double covering of a curve~$\BB$
as above has a type~$\bA_p$ singular point.
\endlemma

\proof
If $p$ were even, the two middle curves in the exceptional divisor
over~$P$ would intersect transversally at a fixed point of~$\tc$.
Since $\tc$ is symplectic, it cannot transpose two such curves.
(The differential $d\tc$ at each fixed point is the multiplication
by~$(-1)$.) Hence, $p$ is odd.

The middle curve in the exceptional divisor is fixed by~$\tc$;
hence, it contains two fixed points. Blowing them up, one obtains
a collection of rational curves with the following incidence
graph:
$$
\DDcenter
@.@.@(\one)\cr
@.@.@|---\cr
@(\two)\withdots\length7@---@(\two)@---@(\four)@---
 @(\two)\withdots\length7@---@(\two)\cr
@.@.@|---\cr
@.@.@(\one)
\endDD
$$
(Here, $\DDcenter@(\one)\endDD$\,, $\DDcenter@(\two)\endDD$\,, and
$\DDcenter@(\four)\endDD$
stand, respectively, for $(-1)$-, $(-2)$-, and $(-4)$-curves; the
total number of curves is $p+2=2k+1$, the action of~$\tc_*$ on
the graph is the horizontal
symmetry, and $\tc$ fixes the two $(-1)$-curves pointwise.)
The projection of this
divisor to~$Y$ is a collection of $(k+2)$ rational curves whose
incidence graph is $\bD_{k+2}$ (the `short' edges corresponding to
the two $(-1)$-curves in~$\tX$).

Since $\tau_Y$ is anti-symplectic, whenever two invariant
curves intersect transversally at a fixed point, exactly one of
them is fixed by~$\tau_Y$ pointwise. Hence, the action of~$\tau_Y$
on the exceptional divisor is determined by its action on the
curve with three neighbors in the incidence graph. Depending on
whether this action is trivial or not, the projection of the
exceptional divisor to $Y\!/\tau_Y$ has one of the following two
incidence graphs:
$$
\DDcenter
@(\one)\cr
@|---\cr
@(\four)@---@(\one)@---@(\four)\length2@---\withdots\length4@---\cr
@|---\cr
@(\one)
\endDD\qquad\text{or}\qquad
\DDcenter
@(\two)@---@(\one)@---@(\four)@---@(\one)\length2@---\withdots\length4@---
\endDD
$$
Now, depending on the parity of~$k$, these graphs blow down either
to a single type~$\bD_{k+2}$ singular point of the branch locus
(if the last vertex corresponds to
a $(-1)$-curve) or to a $(-2)$-curve with a
type~$\bD_{k+1}$ singular point on it (if the last vertex
corresponds to a $(-4)$-curve); in the latter case, the remaining
$(-2)$-curve is the exceptional section $E\subset\Sigma_2$.

The converse statement is proved by analyzing a local equation.
\endproof

For completeness, we also consider the case of type~$\bD$
and~$\bE$ singular points.

\lemma\label{image.D-E}
Let $P$ be a $c$-invariant singular point
of~$B$ of type~$\bD$ or~$\bE$,
%and denote by~$\PP$ its image in~$\Sigma_2$.
and let $\PP\in\Sigma_2$ be its image.
Then $\tc_*$ acts
nontrivially on~$D_P$\rom; in particular, $P$ is of type~$\bD_q$,
$q\ge4$,
or~$\bE_6$. Furthermore, $P\in L$ and~$\PP$ is, respectively, a
type~$\bD_{2q-2}$ or~$\bE_7$ singular point of $\BB+\LL$\rom; in
the former case, $\PP$ is a simple node of~$\BB$ with one of the
branches tangent to~$\LL$.
Conversely, the double covering of a curve~$\BB$ as above has
a corresponding type~$\bD_q$ or~$\bE_6$ singular point.
\endlemma

\proof
If $\tc_*$ acted trivially on~$D_P$, then the curve with three
neighbors in the diagram would have three fixed points and thus
it would be fixed by~$\tc$. Hence, the action is nontrivial. This
observation rules out type~$\bE_7$ and~$\bE_8$ singular points.
The further analysis is completely similar to the proof of
Lemma~\ref{image.nontrivial}, with the additional simplification
that the diagram in~$Y$ is asymmetric and, hence, the curve with
three neighbors is fixed by~$c_Y$.
We omit the details.
\endproof

\subsection{Construction of the involution}\label{involution}
All \term-sextics are described
in~\cite{degt.Oka}; any such curve has `essential' set of
singularities $4\bA_4$, $\bA_9\oplus2\bA_4$, or $2\bA_9$ plus,
possibly, a few other singular points of type~$\bA_1$ or~$\bA_2$.

Let $B$ be a \term-sextic.
Pick a point $P=P_i$ of type~$\bA_4$ (respectively, a point
$P=Q_k$ of
type~$\bA_9$),
%orient
choose an orientation of
its (linear) graph~$D_P$,
and denote by $e_{i1},\ldots,e_{i4}$ (respectively,
$f_{k1},\ldots,f_{k9}$)
%a standard basis for~$\Sigma_P$ formed by
%the classes of the exceptional divisors over~$P$ numbered
%consecutively along the linear graph~$D_P$.
the vertices of~$D_P$, numbered consecutively
%along the graph.
according to the chosen orientation.
Let $e_{ij}^*$
(respectively, $f_{kj}^*$) be the dual basis for
%the dual group
$\Sigma_P^*\subset\Sigma_P\otimes\Q$. Note that
$e_{ij}^*=-e_{i,5-j}^*\bmod\Sigma_P$ and
$f_{ij}^*=-f_{i,10-j}^*\bmod\Sigma_P$.
According
to~\cite{degt.Oka}, under an appropriate numbering of the singular
points and appropriate orientation of their graphs~$D_P$, the
kernel~$\CK$ of the extension $\tSigma\supset\Sigma$ is the cyclic
group~$\CG5$ generated by the
residue $\gamma=\bar\gamma\bmod\Sigma$, where
$\bar\gamma$ is given by
$$
e_{11}^*+e_{21}^*+e_{32}^*+e_{42}^*,\quad
f_{14}^*+e_{11}^*+e_{21}^*,\quad\text{or}\quad
f_{14}^*+f_{22}^*
$$
(for the set of essential singularities
$4\bA_4$, $\bA_9\oplus2\bA_4$, or $2\bA_9$, respectively).

Define an involution $c_S\:S\to S$ as follows:
%$$
%\align
%&h\mapsto h,\\
%&x\mapsto x\quad\text{for $x\in\Sigma_P$ for $P$ a singular point
%other than $\bA_4$ or~$\bA_9$},\\
%&e_{1j}\leftrightarrow e_{2,5-j},\quad
% e_{3j}\leftrightarrow e_{4,5-j},\quad j=1,\ldots,4,\\
%&f_{kj}\leftrightarrow f_{k,10-j},\quad j=1,\ldots,9.
%\endalign
%$$
\Dashes
\roster
\dash
$h\mapsto h$,
\dash
$x\mapsto x$ for $x\in\Sigma_P$ for $P$ a singular point
other than $\bA_4$ or~$\bA_9$,
\dash
$e_{1j}\leftrightarrow e_{2,5-j}$,\quad
$e_{3j}\leftrightarrow e_{4,5-j}$,\quad $j=1,\ldots,4$,
\dash
$f_{kj}\leftrightarrow f_{k,10-j}$,\quad $j=1,\ldots,9$.
\endroster
It is immediate that $c_S$ acts identically on the $p$-primary
part of $\discr S$ for any prime $p\ne5$, and the action of~$c_S$
on the $5$-primary part $\discr S\otimes\Bbb F_5$ (which can be
regarded as an $\Bbb F_5$-vector space) has two dimensional
$(-1)$-eigenspace which contains~$\CK$ as a maximal isotropic
subgroup, so that $\CK^\perp/\CK$ can be identified with the
$(+1)$-eigenspace of~$c_S$.
Hence, $c_S$ extends to an involution $\tc_S\:\tS\to\tS$,
the latter acts identically on $\discr\tS$, and the direct
sum $\tc_S\oplus\id_{S^\perp}$ extends to an involution
$\tc_*\:H_2(X)\to H_2(X)$.

By construction, $\tc_*$ preserves~$h$ and, since it acts
identically on the transcendental lattice
$(\Pic X)^\perp\subset\tS^\perp$, it also preserves classes of
holomorphic forms. Furthermore, $\tc_*$ preserves the positive
cone of~$X$. (Recall that the positive cone is an open
fundamental
polyhedron $V^+\subset(\Pic X)\otimes\R$ of the group generated by
reflections defined by vectors $x\in\Pic X$ with $x^2=-2$;
it is uniquely
characterized by the requirement that $V^+\cdot e>0$ for any
exceptional divisor $e\in\bigcup_PD_P$
and that the closure of~$V^+$ should contain~$h$.)
Due to the description of the fine period space of
marked $K3$-surfaces
given in A.~Beauville~\cite{Beauville},
$\tc_*$ is induced by a unique involutive automorphism
$\tc\:X\to X$,
which is symplectic and commutes with
the deck translation~$\tau$. The descent of~$\tc$
to~$\Cp2$ is the involution~$c$ the existence of which
is asserted by
Theorem~\ref{th.reduction}.

\subsection{Proof of Theorem~\ref{th.reduction}}\label{s.proof}
The involution $c\:\Cp2\to\Cp2$ is constructed in~\ref{involution}.
Due to
Lemmas~\ref{image.class}--\ref{image.nontrivial} (and the
description of the singularities of~$B$ and the action of~$\tc_*$
on the set
%$D=\bigcup_PD_P$
of exceptional divisors,
see~\ref{involution}), the image
$\BB\subset\Sigma_2$ is either a trigonal curve with the set of
singularities
$2\bA_4$ (and then $O\notin B$)
or a hyperelliptic curve with the set of singularities
$2\bA_4$ or $\bA_4\oplus\bA_3$. The latter possibility is ruled
out by the fact that
the genus of a nonsingular curve in $\ls|2E+6F|$ is~$3$.
\qed

\section{Calculation of the groups}\label{S.groups}

\subsection{The curve $\BB$}\label{s.BB}
The trigonal curve $\BB\subset\Sigma_2$ with two type~$\bA_4$
singular points is a maximal trigonal curve in the sense
of~\cite{degt.kplets}. Up to automorphism of~$\Sigma_2$, such a
curve is unique; its skeleton $\Sk\subset\Cp1$
(see~\cite{degt.kplets}) is shown in
Figure~\ref{fig.skeleton}.
One can observe that the skeleton is
symmetric with respect to the dotted grey line (the real structure
$z\mapsto\bar z$ on~$\Cp1$)
and, properly
drown, it is also symmetric with respect to the holomorphic
involution $z\mapsto-1/z$. Hence, the curve~$\BB$ can be chosen
real and symmetric with respect to a real holomorphic involution
of~$\Sigma_2$ (see Section~\ref{S.equations} below for explicit
equations). Furthermore, all singular fibers of~$\BB$ (two cusps
and two vertical tangents) are also real.

\midinsert
\centerline{\picture{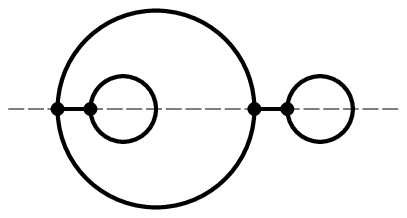}}
\figure\label{fig.skeleton}
The skeleton of~$\BB$
\endfigure
\endinsert

Alternatively, $\BB$ can be obtained as a birational
transform of a plane quartic~$C$ with the set of
singularities $\bA_4\oplus\bA_2$, see Figure~\ref{fig.quartic}.
(Up to automorphism of~$\Cp2$, such a quartic is also unique.)
In the figure, the line $(P_0P_1)$ is tangent to~$C$ at~$P_0$, and
the transformation consists in blowing~$P_0$ up twice and blowing
down the transform of $(P_0P_1)$ and one of the exceptional
divisors over~$P_0$. Lines through~$P_0$ other than $(P_0P_1)$
transform to fibers of~$\Sigma_2$; this observation
gives one a fairly good
understanding of the geometry of~$\BB$, see, \eg,
Figure~\ref{fig.2a9}.

\midinsert
\centerline{\picture{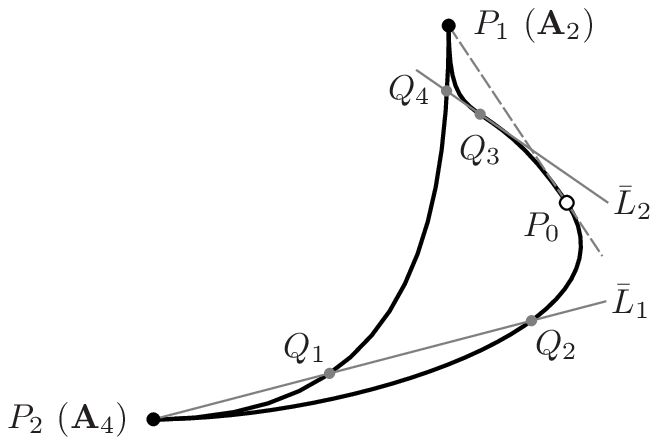}}
\figure\label{fig.quartic}
The quartic~$C$ with the set of singularities $\bA_4\oplus\bA_2$
\endfigure
\endinsert

\subsection{Van Kampen's method}\label{s.vanKampen}
To calculate the fundamental group, we fix a real curve~$\BB$ as
in~\ref{s.BB} and choose an appropriate real section~$\LL$
intersecting~$\BB$ at real points. Let $F_1,\ldots,F_k$ be the
singular fibers of $\BB+\LL$ (\ie, the fibers intersecting
$\BB+\LL$ at less than four points). Under the assumptions, they
are all real. In the figures below, the curves~$\BB$ and~$\LL$ are
shown, respectively, in black and grey, and the singular fibers
are the vertical grey dotted lines.

Fix a real nonsingular fiber~$F_\infty$ intersecting~$\BB$ at one
real point, and consider the affine part
$\Sigma_2\sminus(E\cup F_\infty)$, \cf. Figure~\ref{fig.2a9}. Pick
a real nonsingular fiber~$F$
intersecting~$\BB$ at three real points and
a generic real section~$S$. We identify~$S$ with the base of the
ruling.
Let $x=S\cap F$,
$x_\infty=S\cap F_\infty$, and $x_i=S\cap F_i$, $i=1,\ldots,k$.
Assume that $S$ is \emph{proper} in the following sense: there is
a segment $I\subset S_{\R}$ containing~$x$ and all~$x_i$,
$i=1,\ldots,k$, and disjoint from~$\BB$ and~$\LL$. (The usual
compactness argument shows that such a section exists.) In the
figures, we assume that $I$ lies above~$\BB$ and~$\LL$.

\midinsert
\centerline{\picture{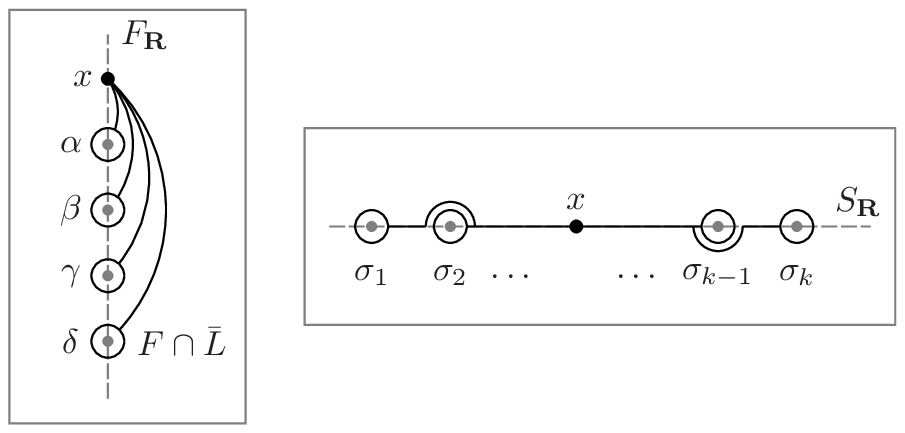}}
\figure\label{fig.basis}
The basis $\Ga$, $\Gb$, $\Gg$, $\Gd$ and the loops~$\Gs_i$
\endfigure
\endinsert

Let $G=\pi_1(F\sminus(\BB\cup E\cup\LL),x)$, and let
$\Ga$, $\Gb$, $\Gg$, $\Gd$ be the basis for~$G$ shown in
Figure~\ref{fig.basis}, left. (All loops are oriented in the
counterclockwise direction.) {\it We always assume
that $\Gd$ is a loop around $F\cap\LL$.}
(Sometimes, the generators should be reordered by inserting~$\Gd$
between $\Ga$ and~$\Gb$ or $\Gb$ and~$\Gg$, \cf.~\ref{s.2a9}
and~\ref{s.4a4+2a1}. In these cases, we show that $\Gd$
commutes with
all subsequent generators, so that the reordering is irrelevant.)
%(In one of the examples,
%see~\ref{s.4a4+2a1}, the correct order is $\Ga$, $\Gb$, $\Gd$,
%$\Gg$. However, in that case $\Gg$ and~$\Gd$ commute and their
%order is irrelevant.)
Let, further, $\Gs_1,\ldots,\Gs_k$ be the
basis for the group
$\pi_1(S\sminus\{x_1,\ldots,x_k\,x_\infty\},x)$ shown in
Figure~\ref{fig.basis}, right: each~$\Gs_i$ is a small circle
about~$x_i$ connected to~$x$ by a real segment~$l_i\subset S_{\R}$
bent to circumvent
the other singular fibers
in the counterclockwise direction.

\definition
The \emph{braid monodromy} along a loop~$\Gs_i$, $i=1,\ldots,k$,
is the automorphism $m_i\:G\to G$ resulting from
dragging~$F$ along~$\Gs_i$ while keeping
the base point in~$S$.
(Since $S$ is proper, $m_i$ is indeed a braid.)
\enddefinition

\proposition\label{th.vanKampen}
The group $\Pi=\pi_1(\Sigma_2\sminus(\BB\cup E\cup\LL))$ is given
by
%the representation
$$
\Pi=\bigl<\Ga,\Gb,\Gg,\Gd\bigm|
 \text{$m_i=\id$, $i=1,\ldots,k$,\quad
 $(\Ga\Gb\Gg\Gd)^2=1$}\bigr>,
 \eqtag\label{eq.vanKampen}
$$
where each \emph{braid relation} $m_i=\id$
%, $i=1,\ldots,k$,
should be understood as a quadruple of relations
$m_i(\Ga)=\Ga$, $m_i(\Gb)=\Gb$, $m_i(\Gg)=\Gg$, $m_i(\Gd)=\Gd$.
\endproposition

\proof
The representation~\eqref{eq.vanKampen} is the essence of van
Kampen's method, see~\cite{vanKampen}, applied to the ruling
of~$\Sigma_2$. The only statement that needs proof is the last
relation $(\Ga\Gb\Gg\Gd)^2=1$, resulting from the patching of the
fiber at infinity~$F_\infty$. This relation is $[\partial D]=1$,
where $D\subset S$ is a small disk around $F_\infty\cap S$. If the
base fiber~$F$ is sufficiently close to~$F_\infty$, then
$[\partial D]=(\Ga\Gb\Gg\Gd)^2$, as in this case one can take
for~$S$ a small perturbation of $E+2F$. In general, due to the
properness of~$S$, the translation homomorphism between any two
nonsingular fibers is a braid; hence, it leaves the product
$\Ga\Gb\Gg\Gd$ invariant and
the relation has the same
form for any fiber~$F$.
\endproof

In sections~\ref{s.2a9}--\ref{s.a9+2a4+a1} below, we attempt to
calculate~$\Pi$ using Proposition~\ref{th.vanKampen}.
To find the braid monodromy~$m_i$, we represent it as the
\emph{local braid monodromy} along a small circle
surrounding~$x_i$, conjugated by the translation homomorphism
along the real path~$l_i$ connecting~$x_i$ to~$x$. The former is
well known: it can be found by considering model
equations. For the latter, we choose the models so that,
at each moment,
all but at most two points of the curve are real; in this case,
the resulting braids are written down directly from the
pictures.

The following lemma facilitates the calculation by reducing the
number of fibers to be considered.

\lemma\label{-1rel}
In representation~\eqref{eq.vanKampen},
\rom(any\rom) one of
the braid relations $m_i=\id$ can be ignored.
\endlemma

\proof
The product $\Gs_1\ldots\Gs_k$ is the class of a
large circle encompassing all singular fibers.
Hence, $m_k\circ{}\ldots{}\circ m_1$ is the so called
\emph{monodromy at infinity}, which is known to be the conjugation
by $(\Ga\Gb\Gg\Gd)^2$. In view of the last relation
in~\eqref{eq.vanKampen},
each~$m_i$ is a composition of the others.
\endproof

For the rest of this section, we fix the notation $\Ga$, $\Gb$,
$\Gg$, $\Gd$ for generators of the group
$\Pi=\pi_1(\Sigma_2\sminus(\BB\cup E\cup\LL))$, chosen as
explained above. The generator~$\Gd$ plays a special r\^ole in the
passage to $\pi=\pi_1(\Cp2\sminus B)$.

\proposition\label{pi}
The fundamental group
$\pi=\pi_1(\Cp2\sminus B)$ is the kernel of the homomorphism
$\Pi\!/\Gd^2\to\CG2$, $\Ga,\Gb,\Gg\mapsto0$, $\Gd\mapsto1$.
\endproposition

\proof
The statement is a direct consequence of the construction: one
considers the appropriate double covering of
$\Sigma_2\sminus(\BB\cup E\cup\LL)$ and patches~$L$.
\endproof

\lemma\label{[]=1}
%In the notation above, let
%$\Pi=\langle\Ga,\Gb,\Gg,\Gd\,|\,\text{relations}\rangle$ be the
%group $\pi_1(\Sigma_2\sminus(\BB\cup E\cup\LL))$, and assume
If
$\Gd$ is a central element of~$\Pi\!/\Gd^2$, then
$\pi_1(\Cp2\sminus B)=\DG{10}\times\CG3$.
\endlemma

\proof
Since $\Gd$ is a central element, the relation
$(\Ga\Gb\Gg\Gd)^2=1$ (or similar) turns into $(\Ga\Gb\Gg)^2=1$ in
$\Pi\!/\Gd^2$, and each braid relation becomes either trivial or a
braid relation for the group
$\pi_1(\Sigma_2\sminus(\BB\cup E))=\DG{10}\times\CG3$
(the latter group is
calculated in~\cite{degt.Oka2}). Hence,
$\Pi\!/\Gd^2=(\DG{10}\times\CG3)\times\CG2$,
%and the statement follows from
and Proposition~\ref{pi} applies.
\endproof

\corollary\label{4a4}
Let~$B$ be a \term-sextic with the set of singularities $4\bA_4$.
Then one has $\pi_1(\Cp2\sminus B)=\DG{10}\times\CG3$.
\endcorollary

\proof
The curve~$B$ is the double covering of~$\BB$
ramified at a section~$\LL$ transversal to~$\BB$. In this case,
$\Gd$ is a central element of
$\pi_1(\Sigma_2\sminus(\BB\cup E\cup\LL))$ (\cf.
Section~\ref{s.2a9} and
Figure~\ref{fig.2a9} below for a much less generic situation), and
the statement follows from Lemma~\ref{[]=1}.
\endproof

\theorem\label{th.converse}
Let $\BB\subset\Sigma_2$ be a trigonal curve with two type~$\bA_4$
singular points, and let $p\:\Cp2\to\Sigma_2/E$ be the double
covering ramified at the vertex $E/E$ and a section~$\LL$ disjoint
from~$E$. Then
%the pull-back
$p^{-1}(\BB)$ is a generalized \term-sextic.
%plane sextic whose
%fundamental group factors to~$\DG{10}$.
\endtheorem

\proof
By perturbing~$\LL$ to a section transversal to~$\BB$, one
perturbs~$B$ to a generic \term-sextic as in Corollary~\ref{4a4}.
\endproof

\subsection{The set of singularities $2\bA_9$}\label{s.2a9}
The sextic~$B$ is the double covering of $\BB$ ramified at a
section~$\LL$ passing through both cusps of~$\BB$, see
Figure~\ref{fig.2a9}; one can take for~$\LL$ the transform of the
secant~$\LL_1$ shown in Figure~\ref{fig.quartic}.

\midinsert
\centerline{\picture{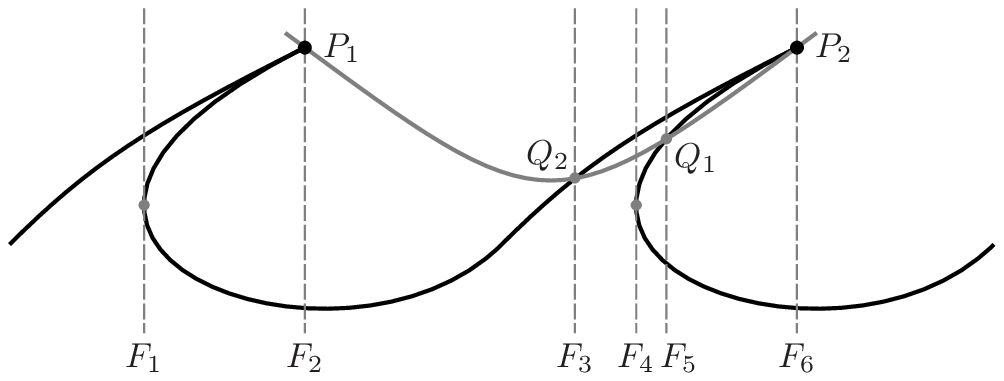}}
\figure\label{fig.2a9}
The set of singularities $2\bA_9$
\endfigure
\endinsert

We choose the generators~$\Ga$, $\Gd$, $\Gb$, $\Gg$ in a
nonsingular fiber between~$F_4$ and~$F_5$. Then, there are
relations $[\Gg,\Gb]=1$, $\Gb=\Gg$, and
$[\Gd,\Ga]=1$ (from~$F_5$, $F_4$, and~$F_3$, respectively); due to
Lemma~\ref{[]=1}, the group is~$\DG{10}\times\CG3$.

\subsection{The set of singularities $4\bA_4\oplus2\bA_1$}\label{s.4a4+2a1}
The curve~$B$ is the double covering of~$\BB$ ramified at a
section~$\LL$ double tangent to~$\BB$,
see Figure~\ref{fig.4a4+2a1}. (The existence of a double tangent
section whose
position with respect to~$\BB$ is as shown in the figure is
rather obvious geometrically: one moves a sufficiently sharp parabola to
achieve two tangency points.
An explicit construction of a pair~$(\BB,\LL)$ using equations
is found in
Section~\ref{s.2tangent} and Figure~\ref{fig.plot} below.)

\midinsert
\centerline{\picture{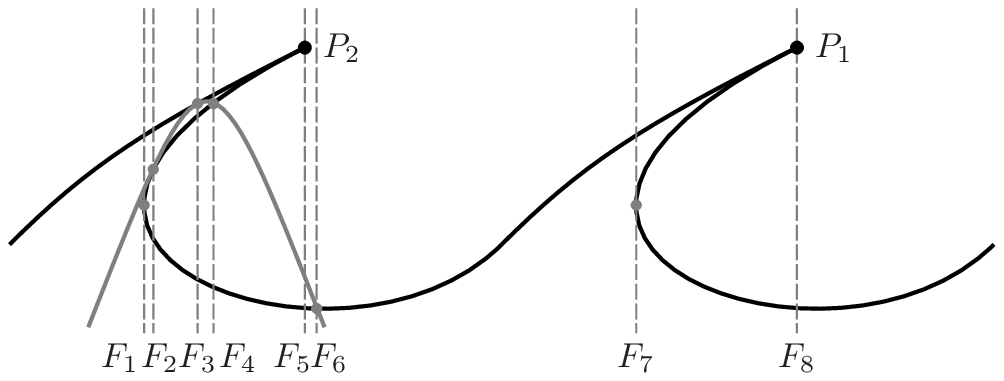}}
\figure\label{fig.4a4+2a1}
The set of singularities $4\bA_4\oplus2\bA_1$
\endfigure
\endinsert

We choose the generators $\Ga$, $\Gb$, $\Gd$, $\Gg$
in a nonsingular fiber between~$F_4$ and~$F_5$.
Ignoring the cusp~$F_8$, see Lemma~\ref{-1rel}, the relations
for~$\Pi$ are
$$\openup\jot
\halign{\thetag{\counter\equation}#&
 \qquad\quad$\displaystyle#$\hfil\qquad&#\hfil\cr
\label{eq2.1}&
 (\Ga\Gb)^2\Ga=\Gb(\Ga\Gb)^2&
 (from the cusp~$F_5$),\cr
\label{eq2.2}&
 [\Gb,\Gd]=[\Gg,\Gd]=1&
 (from~$F_4$ and~$F_6$),\cr
\label{eq2.3}&
 (\Ga\Gd)^2=(\Gd\Ga)^2&
 (from~$F_3$),\cr
\label{eq2.4}&
 (\Ga\1\Gd\Ga\Gb)^2=(\Gb\Ga\1\Gd\Ga)^2&
 (from~$F_2$),\cr
\label{eq2.5}&
 \Gg=(\Ga\1\Gd\Ga)\1\Gb(\Ga\1\Gd\Ga)&
 (from the tangent~$F_1$),\cr
\label{eq2.6}&
 \Gg\1\Ga\Gb\Ga\1\Gg=(\Ga\Gb)\Ga(\Ga\Gb)\1&
 (from the tangent~$F_7$),\cr
\label{eq2.7}&
 (\Ga\Gb\Gd\Gg)^2=1&
 (patching~$F_\infty$).
\crcr}
$$
(The relations are simplified using~\eqref{eq2.2}.) The
corresponding
group~$\pi$ given by Proposition~\ref{pi} was analyzed using the \GAP\
software package. According to \GAP, $\pi$ is an iterated
semi-direct product
$\CG3\times((((\CG2\times\QG)\GAPcolon\CG2)\GAPcolon\CG5)\GAPcolon\CG2)$,
and the abelian factors of its derived series are as stated in
Theorem~\ref{th.main}. (Here, $Q_8$ is the order~$8$ subgroup
$\{\pm1,\pm i,\pm j,\pm k\}\subset\Bbb H$.)

%$$
%\gather
%(\Ga\Gb)^2\Ga=\Gb(\Ga\Gb)^2\qquad
% \text{(from the cusp~$F_5$)},\eqtag\label{eq2.1}\\
%[\Gb,\Gd]=[\Gg,\Gd]=1\qquad
% \text{(from~$F_4$ and~$F_6$)},\eqtag\label{eq2.2}\\
%(\Ga\Gd)^2=(\Gd\Ga)^2\qquad
% \text{(from~$F_3$)},\eqtag\label{eq2.3}\\
%(\Ga\1\Gd\Ga\Gb)^2=(\Gb\Ga\1\Gd\Ga)^2\qquad
% \text{(from~$F_2$)},\eqtag\label{eq2.4}\\
%\Gg=(\Ga\1\Gd\Ga)\1\Gb(\Ga\1\Gd\Ga)\qquad
% \text{(from the tangent~$F_1$)},\eqtag\label{eq2.5}\\
%\Gg\1\Ga\Gb\Ga\1\Gg=(\Ga\Gb)\Ga(\Ga\Gb)\1\qquad
% \text{(from the tangent~$F_7$)},\eqtag\label{eq2.6}\\
%(\Ga\Gb\Gg\Gd)^2=1\qquad
% \text{(patching~$F_\infty$)}.\eqtag\label{eq2.7}
%\endgather
%$$
%(Originally, relation~\eqref{eq2.4} reads
%$(\Gd\1\Ga\1\Gd\Ga\Gd\Gb)^2=(\Gb\Gd\1\Ga\1\Gd\Ga\Gd)^2$; it is
%simplified using $[\Gb,\Gd]=1$. Similarly, \eqref{eq2.5} is
%simplified from its original form using~\eqref{eq2.2}.)
%To discuss the right cusp.\mnote{to do}

\subsection{The set of singularities $\bA_9\oplus2\bA_4\oplus\bA_2$}\label{s.a9+2a4+a2}
The curve~$B$ is the double covering of~$\BB$ ramified at a
section~$\LL$ inflection tangent to~$\BB$,
see Figure~\ref{fig.a9+2a4+a2}. One can take for~$\LL$ the
transform of the tangent~$\LL_2$ shown in
Figure~\ref{fig.quartic}; an explicit construction using equations
is found in Section~\ref{s.inflection} and Figure~\ref{fig.plot}.

\midinsert
\centerline{\picture{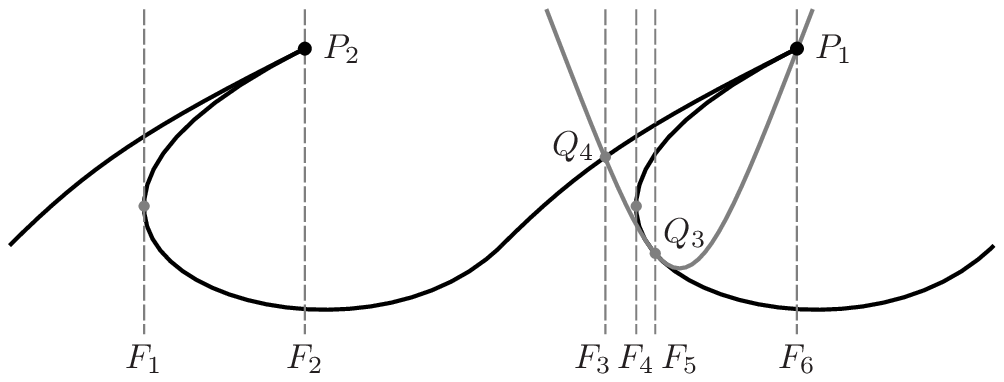}}
\figure\label{fig.a9+2a4+a2}
The set of singularities $\bA_9\oplus2\bA_4\oplus\bA_2$
\endfigure
\endinsert

We choose the generators $\Ga$, $\Gb$, $\Gg$, $\Gd$
in a nonsingular fiber between~$F_4$ and~$F_5$. Ignoring the
vertical tangent~$F_1$, see Lemma~\ref{-1rel}, the relations
for~$\Pi$ are
$$\openup\jot
\halign{\thetag{\counter\equation}#\hss&
 \qquad\quad$\displaystyle#$\hfil\qquad&#\hfill\cr
\label{eq3.1}&
 \Gb=\Gg&
 (from the tangent~$F_4$),\cr
\label{eq3.2}&
 [\Ga,\Gg\Gd\Gg\1]=1&
 (from~$F_3$),\cr
\label{eq3.3}&
 (\Gg\Gd)^3=(\Gd\Gg)^3&
 (from~$F_5$),\cr
\label{eq3.4}&
 [\Ga\Gb,\Gd_1]=1&
 (from the cusp~$F_6$),\cr
\label{eq3.5}&
 \Gd_1(\Ga\Gb)^2\Ga=\Gb(\Ga\Gb)^2\Gd_1&
 (from the cusp~$F_6$),\cr
\label{eq3.6}&
 (\Gb_2\Gg)^2\Gb_2=\Gg(\Gb_2\Gg)^2&
 (from the cusp~$F_2$),\cr
\label{eq3.7}&
 (\Ga\Gb\Gg\Gd)^2=1&
 (patching~$F_\infty$),\cr
\crcr}
$$
%$$
%\gather
%\Gb=\Gg\qquad
% \text{(from the tangent~$F_4$)},\eqtag\label{eq3.1}\\
%[\Ga,\Gg\Gd\Gg\1]=1\qquad
% \text{(from~$F_3$)},\eqtag\label{eq3.2}\\
%(\Gg\Gd)^3=(\Gd\Gg)^3\qquad
% \text{(from~$F_5$)},\eqtag\label{eq3.3}\\
%[\Ga\Gb,\Gd_1]=1\qquad
% \text{(from the cusp~$F_6$)},\eqtag\label{eq3.4}\\
%\Gd_1(\Ga\Gb)^2\Ga=\Gb(\Ga\Gb)^2\Gd_1\qquad
% \text{(from the cusp~$F_6$)},\eqtag\label{eq3.5}\\
%(\Gb_2\Gg)^2\Gb_2=\Gg(\Gb_2\Gg)^2\qquad
% \text{(from the cusp~$F_2$)},\eqtag\label{eq3.6}\\
%(\Ga\Gb\Gg\Gd)^2=1\qquad
% \text{(patching~$F_\infty$)},\eqtag\label{eq3.7}
%\endgather
%$$
where $\Gd_1=(\Gg\Gd\Gg)\Gd(\Gg\Gd\Gg)\1$ and
$\Gb_2=(\Ga\Gg\Gd\1\Gg\1)\Gb(\Ga\Gg\Gd\1\Gg\1)\1$. The resulting
group~$\pi$, see Proposition~\ref{pi},
was analyzed using \GAP. Its derived series is as
stated in Theorem~\ref{th.main}.

\subsection{The sets of singularities $\bA_9\oplus2\bA_4\oplus\bA_1$
and $4\bA_4\oplus\bA_2$}\label{s.a9+2a4+a1}
For the set of singularities $\bA_9\oplus2\bA_4\oplus\bA_1$, we
perturb the inflection tangency point~$Q_3$ in
Figure~\ref{fig.a9+2a4+a2} to a simple tangency point and a point
of transversal intersection. Then, \eqref{eq3.3} is replaced with
$[\Gg,\Gd]=1$ and one obtains $[\Ga,\Gd]=1$ (from~\eqref{eq3.2}\,),
$\Gd_1=\Gd$, and $[\Gb,\Gd]=1$ (from~\eqref{eq3.4}\,); due to
Lemma~\ref{[]=1}, the resulting group~$\pi$ is~$\DG{10}\times\CG3$.

For the set of singularities $4\bA_4\oplus\bA_2$, the intersection
point~$P_1$ in Figure~\ref{fig.a9+2a4+a2} is perturbed to two
points of transversal intersection. Then, \eqref{eq3.4}
and~\eqref{eq3.5} turn into
$[\Ga,\Gd_1]=[\Gb,\Gd_1]=1$ and
$(\Ga\Gb)^2\Ga=\Gb(\Ga\Gb)^2$, respectively. Using \GAP\ shows
that the resulting group~$\pi$ is~$\DG{10}\times\CG3$.

\subsection{The sets of singularities $\bA_9\oplus2\bA_4$ and
$4\bA_4\oplus\bA_1$}\label{s.others}
These sextics can be obtained by small perturbations from sextics
with the sets of singularities, \eg, $2\bA_9$
(see~\ref{s.2a9})
and
$\bA_9\oplus2\bA_4\oplus\bA_1$ (see~\ref{s.a9+2a4+a1}),
respectively; the resulting groups
have already been shown to be isomorphic to $\DG{10}\times\CG3$.

\section{The equations}\label{S.equations}

The calculations in this section (substitution, factorization,
discriminants, and system solving) are straightforward but rather
tedious. Most calculations were performed using \Maple.

\subsection{The curve~$\BB$}
In appropriate affine coordinates $(x,y)$ in~$\Sigma_2$ the
trigonal curve~$\BB$ with two type~$\bA_4$ singular points is
given by the Weierstra{\ss} equation
$$
\gather
f(x,y)=4y^3-3yp(x)+q(x)=0,\ \text{where}\eqtag\label{eq.B}\\
p(x)=x^4-12x^3+14x^2+12x+1,\\
q(x)=(x^2+1)(x^4-18x^3+74x^2+18x+1).
\endgather
%4y^3-3y(x^4-12x^3+14x^2+12x+1)+(x^2+1)(x^4-18x^3+74x^2+18x+1)
$$
%where $p=x^4-12x^3+14x^2+12x+1$ and
%$q=(x^2+1)(x^4-18x^3+74x^2+18x+1)$.
The discriminant of~\eqref{eq.B} with respect to~$y$ is
$\Delta=(2)^{10}(3)^6x^5(x^2-11x-1)$; it has two $5$-fold
roots $x=0$ and $x=\infty$ (the singular points
of~$\BB$)
and two simple roots $x_\pm=11/2\pm5\sqrt5/2$ (the two vertical tangents).

\remark{Remark}
The point $x=\infty$ is a $5$-fold root of~$\Delta$
as the `predicted' degree
of~$\Delta$ is~$12$. Originally, equation~\eqref{eq.B} was
obtained by an appropriate birational coordinate change from the
equation $y^2-2x^2y+x^4-4x^3y$ of the quartic with the set of
singularities $\bA_4\oplus\bA_2$, see Figure~\ref{fig.quartic}.
\endremark

The curve~$\BB$ is rational; it can be parametrized as follows
$$
x(t)=\frac{t^2(t-1)}{t+1}\,,\qquad
 y(t)=\frac{(t^2+1)(t^4-2t^3-6t^2+2t+1)}{2(t+1)^2}\,.
\eqtag\label{eq.parametric}
$$
The special points on the curve correspond to the following values
of the parameter:
$$
\alignedat2
&t_0=0,\quad t_\infty=\infty&\qquad&
 \text{(the cusps)},\\
&t_0'=1,\quad t_\infty'=-1&&
 \text{(the points under the cusps)},\\
&t_\pm=-\frac12\mp\frac12\sqrt5&&
 \text{(the tangency points over $x=x_\pm$)},\\
&t_\pm'=2\pm\sqrt5&&
 \text{(the other points over $x=x_\pm$)}.
\endalignedat\eqtag\label{eq.t}
$$

Both the curve and the parametrization are real, as are all
singular fibers of~$\BB$. Furthermore, $\BB$ is invariant under
the automorphism $x\mapsto-1/x$, $y\mapsto y/x^2$. In the
$t$-line, this transformation corresponds to the automorphism
$t\mapsto-1/t$.

\subsection{Generic sextics}
Due to Theorem~\ref{th.reduction}, any \term-sextic
is given by an affine equation of the form
$$
f(x, y^2+x^2+bx+c)=0,\eqtag\label{eq.sextic}
$$
where $f(x,y)$ is the polynomial given by~\ref{eq.B} and
$y=ax^2+bx+c$ is the equation of the section~$\LL$ constituting
the branch locus. Conversely, from Theorem~\ref{th.converse} it
follows that any curve~$B$ given by~\eqref{eq.sextic} is a
\term-sextic provided that it is irreducible and all its
singularities are simple. Note that $B$ is reducible (splits into
two cubics interchanged by the involution on~$\Cp2$) if and only
if, at each point of intersection of~$\BB$ and~$\LL$, the local
intersection index is even. Hence, in view of the classification
of sections given below, $B$ is reducible if and only if it has
the (non-simple) set of singularities $\bY_{1,1}^1\oplus\bA_9$,
see~\ref{s.sing}; such a curve splits into two cubics with a
common cusp.

The set of singularities of a sextic~$B$ given
by~\eqref{eq.sextic} with a generic triple $(a,b,c)$ (so that
$\LL$ is transversal to~$\BB$) is $4\bA_4$. In
Sections~\ref{s.sing}--\ref{s.inflection} below, we discuss the
possible degenerations of the section~$\LL$ and express them in
terms of the triple $(a,b,c)$. (Sometimes, the condition is stated
using an extra parameter~$t$,
as an attempt to eliminate~$t$ results in a
multi-line \Maple\ output.) For each degeneration, we use
Lemmas~\ref{image.trivial}
and~\ref{image.nontrivial} to indicate the set of singularities of
the corresponding sextic~$B$. The results should be understood as
follows: a sextic $B$ given by~\eqref{eq.sextic} has a certain set
of singularities~$\Sigma$ if and only if the triple $(a,b,c)$
satisfies the condition corresponding to~$\Sigma$ but does not
satisfy any condition corresponding to an immediate
degeneration of~$\Sigma$
(see the adjacency diagram shown in Figure~\ref{fig.adj}).

\midinsert
$$
\eightpoint\minCDarrowwidth{.65cm}
\CD
@.@.(\bW_{12}\oplus2\bA_4)\\
@.@.@VVV\\
@.@.(\bY_{1,1}^1\oplus2\bA_4)@<<<(\bY_{1,1}^1\oplus\bA_9)\\
@.@.@VVV@VVV\\
@.4\bA_4@<<<\bA_9\oplus2\bA_4@<<<2\bA_9\\
@.@AAA@AAA\\
4\bA_4\oplus2\bA_1@>>>4\bA_4\oplus\bA_1@<<<\bA_9\oplus2\bA_4\oplus\bA_1\\
@.@AAA@AAA\\
@.4\bA_4\oplus\bA_2@<<<\bA_9\oplus2\bA_4\oplus\bA_2\\
\endCD
$$
\figure\label{fig.adj}
Immediate adjacencies of sets of singularities
\endfigure
\endinsert

For completeness, we also mention (parenthetically in
Figure~\ref{fig.adj}) the sextics~$B$ given by~\eqref{eq.sextic} whose
singularities are not simple; this is the case if and only if the
triple $(a,b,c)$ is as in~\eqref{eq.nonsimple} below. In
Arnol$'\!$d's notation, $B$ may only have a non-simple singular
point of one of the following two types:
$\bY^1_{1,1}$ (transversal
intersection of two cusps) or $\bW_{12}$ (semiquasihomogeneous
singularity of type $(4,5)$\,).

\subsection{Sections through singular points}\label{s.sing}
A section $y=ax^2+bx+c$ passes through one of the cusps of~$\BB$
(the set of singularities $\bA_9\oplus2\bA_4$) if and only if
$$
c=\frac12\quad\text{(the cusp at $t=0$)}\qquad\text{or}\qquad
a=\frac12\quad\text{(the cusp at $t=\infty$)}.
\eqtag\label{eq.cusp}
$$
Hence, the section passes through both cusps (the set of
singularities $2\bA_9$) if and only if $a=c=1/2$.

Further degenerations considered here produce
sextics with a non-simple singular point. The section is tangent
to~$\BB$ at a cusp (the set of singularities
$\bY_{1,1}^1\oplus2\bA_4$) if and only if
$$
c=\frac12,\ b=3\quad\text{(at $t=0$)}\qquad\text{or}\qquad
a=\frac12,\ b=-3\quad\text{(at $t=\infty$)}.
\eqtag\label{eq.nonsimple}
$$
It passes through the other cusp (the set of singularities
$\bY_{1,1}^1\oplus\bA_9$; this sextic is reducible) if and only if
$a=c=1/2$ and $b=\pm3$. Finally, the section passes through a cusp
with local intersection index~$5$ (the set of singularities
$\bW_{12}\oplus2\bA_4$) if and only if
$$
(a,b,c)=\Bigl(-\frac{11}2, 3, \frac12\Bigr)
\quad\text{or}\quad
(a,b,c)=\Bigl(\frac12, -3, -\frac{11}2\Bigr).
\eqtag\label{eq.W}
$$
Such a section cannot pass through the other cusp.

A section passing through both cusps of~$\BB$ or tangent to~$\BB$
at a cusp does not admit any degenerations other than described
above. Indeed, if $a=c=1/2$ (respectively, $c=1/2$ and $b=3$),
then, restricting the original polynomial $f(x,y)$ to the section
and reducing the trivial factor~$x^2$ (respectively,~$x^4$), one
obtains a polynomial in~$x$ whose discriminant is
$16(b-3)^3(b+3)^3$ (respectively, $12(2a-1)^5$).

\remark{Remark}
According to~\cite{degt.Oka}, \term-sextics are characterized by
the existence of two conics in a very special position with
respect to the type~$\bA_4$ and~$\bA_9$ singular points of the
curve. These conics are the pull-backs of the two sections
$y=ax^2+bx+c$ with $a=c=1/2$ and $b=\pm3$, each section being
tangent to~$\BB$ at one of its cusps
and passing through the other cusp.
\endremark

\subsection{Digression: other
generalized \term-sextics}\label{s.non-simple}
From~\ref{s.sing},
it follows that the double covering construction also produces
representatives of the two families of irreducible
generalized \term-sextics
with a quadruple singular point, see~\cite{degt.Oka}. (In each
family, symmetric curves form a codimension one subset.)
It is worth mentioning that the remaining classes of irreducible
generalized \term-sextics, those with the sets of singularities
$\bJ_{2,0}\oplus2\bA_4$,
$\bJ_{2,1}\oplus2\bA_4$, and $\bJ_{2,5}\oplus\bA_4$,
see~\cite{degt.Oka2}, are also related to the trigonal curve
$\BB\subset\Sigma_2$ with two type~$\bA_4$ singular points: they
are obtained from~$\BB$ by a birational transformation rather than
double covering.

\subsection{Simple and double tangents}\label{s.tangent}
Let
$$
s(t)=ax^2(t)+bx(t)+c,\eqtag\label{eq.section}
$$
where $x(t)$ is given by~\eqref{eq.t}.
Solving $s(t)=y(t)$ and $s'(t)=y'(t)$, one concludes that
a section $y=ax^2+bx+c$ is tangent to~$\BB$ at a point
$(x(t), y(t))\in\BB$ (the set of singularities
$4\bA_4\oplus\bA_1$) if and only if
$$
\aligned
a&=\frac{-b(t^3+2t^2-1)+t^5-5t^3-5t^2-3}{2t^2(t-1)(t^2+t-1)},\\
c&=-\frac{bt^2(t^3-2t+1)+3t^5+5t^3-5t^2+1}{2(t+1)(t^2+t-1)}
\endaligned\eqtag\label{eq.tangent}
$$
for some $b\in\C$ and $t\in\C\sminus\{0,\pm1,t_\pm\}$ or
$$
\text{$t=1$ and $(b,c)=(-6,-1)$\qquad or\qquad
$t=-1$ and $(a,b)=(-1,6)$.}\eqtag\label{eq.tangent.ex}
$$
This section passes through the cusp of~$\BB$ at $t=0$ (the set of
singularities $\bA_9\oplus2\bA_4\oplus\bA_1$)
if and only if $c=1/2$, see~\ref{s.sing}; in this case
$$
a=\frac{t^4-t^3-2t^2+3t+11}{2(t-1)^2(t^2+t-1)},\quad
b=-\frac{3(t^3+2t-1)}{(t-1)(t^2+t-1)},\quad
c=\frac12.\eqtag\label{eq.tangent0}
$$
The section passes through the cusp of~$\BB$ at $t=\infty$
(another implementation of the set of
singularities $\bA_9\oplus2\bA_4\oplus\bA_1$)
if and only if
$a=1/2$, see~\ref{s.sing}; in this case
$$
a=\frac12,\quad
b=-\frac{3(t^3+2t^2+1)}{(t+1)(t^2+t-1)},\quad
c=-\frac{11t^4-3t^3-2t^2+t+1}{2(t+1)^2(t^2+t-1)}.\eqtag\label{eq.tangent1}
$$
Relations~\eqref{eq.tangent0} and~\eqref{eq.tangent1} still hold
for the exceptional values $t=-1$ and $t=1$, respectively,
\cf.~\eqref{eq.tangent.ex}.

There is a somewhat unexpectedly simple relation between the two tangency
points of a section double tangent to~$\BB$. We state it below as a
separate lemma. Denote by~$\Ge_\pm$ the roots of the polynomial
$t^2+3t+1$. One has
$\Ge_\pm=(-3\pm\sqrt5)/2=t_\pm/t'_\pm$. Note
that $\Ge_+\Ge_-=1$.

\lemma
Let $\BB\subset\Sigma_2$ be the trigonal curve parametrized
by~\eqref{eq.parametric}, and let $t_1,t_2\in\C\sminus\{0,t_\pm\}$
be two distinct values of the parameter. Then, there is a section
tangent to~$\BB$ at both points $(x(t_i),y(t_i))$, $i=1,2$, if and
only if $t_2/t_1=\Ge_\pm$.
\endlemma

\proof
Assume that $t_1,t_2\ne\pm1$. (The case when one of~$t_1$, $t_2$
takes an exceptional value~$\pm1$ is treated similarly
using~\eqref{eq.tangent.ex}.) Let $a_1$, $a_2$ and $c_1$, $c_2$ be
the coefficients~$a$ and~$c$ in~\eqref{eq.tangent} evaluated at
$t=t_1, t_2$, respectively. Then $a_1-a_2=c_1-c_2=0$; hence,
$(a_1-a_2)t_1^2t_2^2(t_1-1)(t_2-1)+(c_1-c_2)(t_1+1)(t_2+1)=0$.
The latter expression takes the form
$$
%(a_1-a_2)t_1^2t_2^2(t_1-1)(t_2-1)+(c_1-c_2)(t_1+1)(t_2+1)=
\frac{3(t_1-t_2)^3(t_1^2+3t_1t_2+t_2^2)}
 {(t_1^2+t_1-1)(t_2^2+t_2-1)}=0,
$$
and, taking into account the restrictions on~$t_1$, $t_2$, one
obtains $t_2/t_1=\Ge_\pm$. For the converse statement, one
observes that, if $t_2/t_1=\Ge_\pm$, then the linear system
$a_1=a_2$, $c_1=c_2$ in one variable~$b$ has a solution
(given by~\eqref{eq.2tangent} below).
\endproof

Thus, a section $y=ax^2+bx+c$ is double tangent to~$\BB$ (the set of
singularities $4\bA_4\oplus2\bA_1$)
if and only if
$$
b=b_\pm=-\frac{3(t^2+(3\Ge_\mp+1)t-\Ge_\mp)(t^2-\Ge_\pm t-\Ge_\mp)(t+t_\pm)}
{(t-t_\pm)(t-t_\mp)^2(t-t'_\pm)^2}.\eqtag\label{eq.2tangent}
$$
%%%% Maple output
%$$
%b_+=
%-\frac{3(2t^2-7t-3\sqrt5t+3+\sqrt5)(2t^2+3t-\sqrt5t+3+\sqrt5)(2t-1-\sqrt5)}
%{(2t+1+\sqrt5)(2t+1-\sqrt5)^2(t-\sqrt5-2)^2}
%$$
%$$
%b_-=
%-\frac{3(2t^2-7t+3\sqrt5t+3-\sqrt5)(2t^2+3t+\sqrt5t+3-\sqrt5)(2t-1+\sqrt5)}
%{(2t+1-\sqrt5)(2t+1+\sqrt5)^2(t-2+\sqrt5)^2}
%$$
%%%%
Here, the two tangency points are at~$t$ and $\Ge_\pm t$;
the expressions for $a$ and~$c$ are obtained by a direct
substitution to~\eqref{eq.tangent}. This relation still holds if
$t=\pm1$.
%The relation $t_2/t_1=\Ga_\pm$ and
%expression~\eqref{eq.2tangent} for~$b$ still hold if one of~$t_1$,
%$t_2$ takes an exceptional value~$\pm1$,
%\cf.~\eqref{eq.tangent.ex}.

A section tangent to~$\BB$ at a smooth point does not admit
any other degenerations. Indeed, sections passing through both
singular points of~$\BB$ or tangent to~$\BB$ at a singular point
are considered in~\ref{s.sing}, and sections inflection tangent
to~$\BB$ are treated in~\ref{s.inflection} below. A section cannot
be tangent to~$\BB$ at three smooth points at $t=t_1$, $t_2$,
and~$t_3$, as then one would have $t_2/t_1=\Ge_\pm$,
$t_3/t_2=\Ge_\pm$, and $t_3/t_1=\Ge_\pm$; this system is
incompatible unless $t_1=t_2=t_3=0$. Finally, a double tangent
cannot pass through a singular point, say at $t=0$, as
substituting $t=t_1$ and $t=t_2$ to~\eqref{eq.tangent0} and
eliminating~$a$ and~$b$,
one obtains a system in~$(t_1,t_2)$ which has no solutions with
$t_1\ne t_2$.

\subsection{Digression: the curve in Figure~\ref{fig.4a4+2a1}}\label{s.2tangent}
The pair $(\BB,\LL)$ used to calculate the fundamental group of a
sextic~$B$ with the set of singularities $2\bA_4\oplus2\bA_1$,
see~\ref{s.4a4+2a1} and Figure~\ref{fig.4a4+2a1}, can be obtained
from~\eqref{eq.2tangent} and~\eqref{eq.tangent} with the following
values of the parameters: $t_1=t=1/2$ and
$t_2=\Ge_+t_1\approx-0.191$. Then one has
$a\approx-161.05$, $b\approx-13.93$, and
$c\approx0.0448$. The two points of transversal intersection
of~$\BB$ and the section $y=ax^2+bx+c$ are at
$t\approx0.281$ (over $x\approx-0.0442$) and $t\approx1.101$
(over $x\approx0.0585$). This section is shown
in Figure~\ref{fig.plot}.

\midinsert
\centerline{\epsfxsize.8\hsize\picture{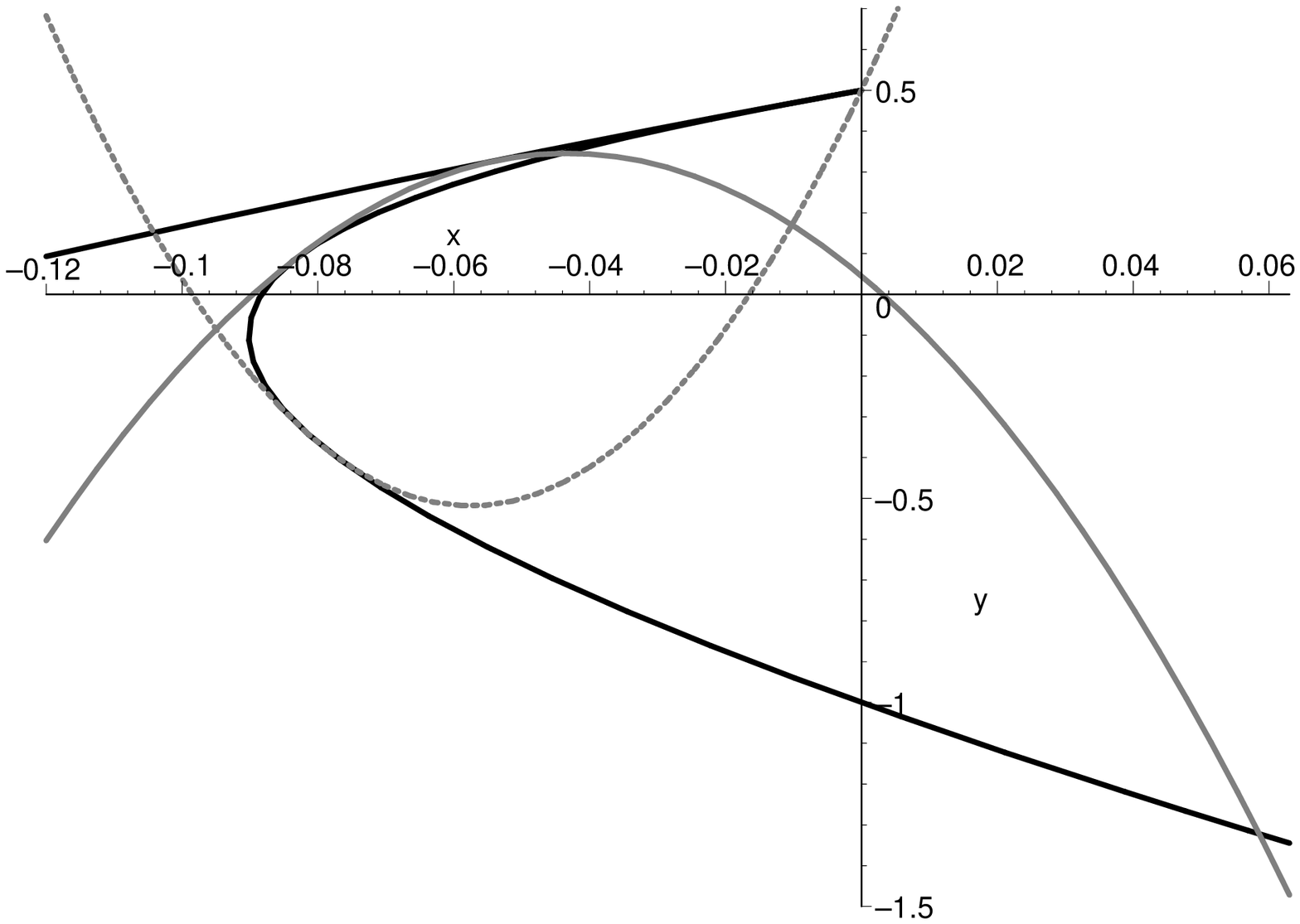}}
\figure\label{fig.plot}
\Maple\ plot of the curve~$\BB$ (black), a double tangent (solid
grey), and an inflection tangent through a singular point (dotted
grey)
\endfigure
\endinsert

\subsection{Inflection tangents}\label{s.inflection}
A section $y=ax^2+bx+c$ is inflection tangent to~$\BB$
at a point $(x(t), y(t))\in\BB$ (the set of singularities
$4\bA_4\oplus\bA_2$) if and only if
$$
\aligned
a&=\frac{t^6+3t^5-5t^3+12t+11}{2(t^2+t-1)^3},\\
b&=-\frac{3(t^2+1)(t^4+3t^3-t^2-3t+1)}{(t^2+t-1)^3},\\
c&=-\frac{11t^6-12t^5+5t^3-3t+1}{2(t^2+t-1)^3},
\endaligned\eqtag\label{eq.infl}
$$
$t\in\C\sminus\{0,t_\pm\}$.
(To see this, one should solve for~$(a,b,c)$ the system
$s(t)=y(t)$, $s'(t)=y'(t)$, $s''(t)=y''(t)$, where $s(t)$ is the
section given by~\eqref{eq.section}.)
This inflection tangent passes through one of the cusps of~$\BB$
(the set of singularities $\bA_9\oplus2\bA_4\oplus\bA_2$) if and
only if $t=3/4$ (the cusp at $t=0$) or $t=-4/3$
(the cusp at $t=\infty$). The corresponding values of $(a,b,c)$
are
$$
(a,b,c)=\Bigl(\frac{3077}{10}, \frac{177}{5}, \frac12\Bigr)
\quad\text{and}\quad
(a,b,c)=\Bigl(\frac12, -\frac{177}5, \frac{3077}{10}\Bigr),
\eqtag\label{eq.infl+cusp}
$$
respectively. The section corresponding to $t=3/4$ is shown in
Figure~\ref{fig.plot}.

An inflection tangent at a smooth point of~$\BB$ cannot have any
other degenerations. Indeed, after clearing the denominators and
reducing the trivial factor $(u-t)^3$, the equation
$y(u)=ax(u)^2+bx(u)+c$ with $a$, $b$, and~$c$ given
by~\eqref{eq.infl} and $x(\,\cdot\,)$, $y(\,\cdot\,)$ as
in~\eqref{eq.parametric}
has solution $u=t$ only for $t=0$ or~$t_\pm$
(hence, no quadruple intersection points),
and the discriminant of the above equation
with respect to~$u$ is, up to a constant
coefficient, $t^2(3t+4)(4t-3)(t^2+t-1)^3$ (hence, no other
tangency points).

\widestnumber\key{EO1}
\refstyle{C}

\enddocument